\documentclass[12pt,oneside]{amsart}
\usepackage[utf8]{inputenc}
\usepackage{geometry}  % See geometry.pdf to learn the layout options. There are lots.
\usepackage{graphicx}
\usepackage{amssymb,amsmath}
\usepackage{epstopdf}
\usepackage[labelformat=empty]{subfig}
\usepackage[pdftex]{color}
\usepackage{xcolor}
\usepackage{bm}
\usepackage{tikz}
\usepackage{url}
\usetikzlibrary{matrix}
\usepackage{multicol}
\usepackage{array}% http://ctan.org/pkg/array
\allowdisplaybreaks
\title{Survey on recent developments in semitoric systems}
\author{Jaume Alonso \qquad Sonja Hohloch}
\numberwithin{equation}{section}
\newtheorem{de}{Definition}[section]

\newtheorem{ex}[de]{Example}

\usepackage[parfill]{parskip}

\newcommand{\vungoc}{V\~u Ng\d{o}c }
\newcommand{\nff}{{n_{\text{FF}}}}

\newcommand{\C}{\mathbb{C}}

\newcommand{\N}{\mathbb{N}}
\newcommand{\R}{\mathbb{R}}
\newcommand{\mbS}{\mathbb{S}}
\newcommand{\T}{\mathbb{T}} % torus
\newcommand{\Z}{\mathbb{Z}}

\newcommand{\al}{\alpha}
\newcommand{\be}{\beta}
\newcommand{\ga}{\gamma}

\newcommand{\eps}{\varepsilon}

\newcommand{\lam}{\lambda}

\newcommand{\om}{\omega}

\newcommand{\Lam}{\Lambda}

\newcommand{\mcC}{\mathcal C}

\newcommand{\mcG}{\mathcal G}

\newcommand{\mcI}{\mathcal I}

\newcommand{\mcL}{\mathcal L}

\newcommand{\mcO}{\mathcal O}

\newcommand{\mcX}{\mathcal X}

\begin{document}

\begin{abstract}
Semitoric systems are a special class of four-dimensional completely integrable systems where one of the first integrals generates an $\mbS^1$-action. They were classified by Pelayo \& \vungoc in terms of five symplectic invariants about a decade ago. We give a survey over the recent progress which has been mostly focused on the explicit computation of the symplectic invariants for families of semitoric systems depending on several parameters and the generation of new examples with certain properties, such as a specific number of singularities of lowest rank.
\end{abstract}

\maketitle
\thispagestyle{empty}
\section{Introduction and motivation}

Semitoric systems are a special class of four-dimensional completely integrable systems where one of the first integrals is proper and generates a global $\mbS^1$-action, singularieties are assumed to be nondegenerate, and none of the singularities has hyperbolic components. These systems can be seen as a generalisation of toric systems with the important difference that only one of the first integrals is required to be periodic. This allows for the existence of focus-focus fibres, which in turn generate monodromy, cf.\ Duistermaat \cite{Du}. The two older surveys by Pelayo \& \vungoc \cite{PV2, PV5} give a good overview of the symplectic classification by Pelayo \& \vungoc \cite{PV1,PV4} and the main features of semitoric systems. Note that semitoric systems appear naturally in physics, for example in the Jaynes-Cummings model (see Babelon \& Cantini \& Douçot \cite{BCD}) and the coupled angular momenta (see Sadovskii \& Zhilinskii \cite{SZ}). However, they have also raised interest in the symplectic-dynamical community since the classification of these systems is in terms of five {\em symplectic} invariants.

This survey aims at giving an overview over the research programme in semitoric systems completed in recent years. This includes the computation of all invariants of the most accessible semitoric system, the coupled spin-oscillator (cf.\ Alonso \& Dullin \& Hohloch \cite{ADH}), the computation of the invariants of the coupled angular momenta, a family of examples depending on three parameters (cf.\ Alonso \& Dullin \& Hohloch \cite{ADH2}), the construction of an example with two focus-focus points (cf.\ Hohloch \& Palmer \cite{HP}) and the generation of examples in all Hirzebruch surfaces (cf.\ Le Floch \& Palmer \cite{LFPa}).

The organisation of the survey is as follows. In the second section, we give the precise definition of semitoric systems and state the classification in terms of the five symplectic invariants. In the third section, we discuss the coupled angular momenta and its symplectic invariants. In the fourth section, we introduce the new examples of semitoric systems, one with two focus-focus points and others in Hirzebruch surfaces. In the fifth section, we give an overview on related research.

\subsection*{Acknowledgements}

The first author has been fully and the second author has been partially funded by the project BOF-DocPro4 FFB150339 of the Research Fund of the University of Antwerp.

\section{The classification and invariants}
In this section, we briefly fix some notations and conventions and recall the definition of a semitoric system and their classification in terms of five symplectic invariants.

\subsection{Integrable systems and semitoric systems}
\label{ssintsys}
Let $(M,\om)$ be a symplectic manifold of dimension $2n$. By making use of the non-degeneracy of $\om$, we associate to any function $f\in \mcC^\infty(M)$ its \emph{Hamiltonian vector field} $\mcX_f$, uniquely defined by the relation $\om(\mcX_f,\cdot) =- df(\cdot)$. The vector field $\mcX_f$ generates a flow on $M$, called \emph{Hamiltonian flow} of $f$. The \emph{Poisson bracket} of two Hamiltonian functions $f,g \in \mcC^\infty(M)$ is defined as $\{f,g\}:=\om(\mcX_f,\mcX_g)$. If $\{f,g\}=0$, we say that $f,g$ \emph{Poisson-commute}, i.e.\ the level sets of each of them are invariant under the flow of the other. A triplet $(M,\om,F)$, where $F=(f_1,...,f_n) \in \C^\infty(M,\R^n)$, is an \emph{integrable system} if $\{f_i,f_j\}=0$ for all $i,j=1,...,n$ and if $\mcX_{f_1},...,\mcX_{f_n}$ are almost everywhere linearly independent. The map $F$, often called \emph{momentum map}, induces a --possibly singular-- fibration on $M$. A point $m\in M$ is called \emph{regular} if $dF$ has maximal rank and \emph{singular} otherwise. The \emph{Arnold-Liouville theorem} fully describes the dynamics on \emph{regular fibres}, the so-called \emph{Liouville tori}. The interesting part, therefore, is to understand the role of the \emph{singular fibres}, i.e.\ those containing at least one singular point. 

We restrict ourselves to \emph{non-degenerate} singular points, i.e.\ points $m\in M$ for which the Hessians $D^2f_1|_m,...,D^2f_n|_m$ span a Cartan subalgebra of the Lie algebra of quadratic forms on $T_mM$ (cf.\ Bolsinov-Fomenko \cite{BF}). These singular points are described by Eliasson \cite{El1,El2} and Miranda \& Zung \cite{MZ} using normal forms. More precisely, for each non-degenerate singularity $m \in M$ we can find an open neighbourhood $m \in U \subseteq M$, symplectic coordinates $(q,p)=(q_1,...,q_n,p_1,...,p_n)$ and a function $G=(G_1,...,G_n):U \to \R^n$ such that $m=(0,0)$, $\{f_i,G_j\}=0$ for all $i,j=1,...,n$ and the components $G_j$ are of the following forms:
\begin{itemize}
	\item Regular component: $G_j(q,p) = p_j$
	\item Elliptic component: $G_j(q,p) = \frac{1}{2}({q_j}^2+{p_j}^2)$
	\item Hyperbolic component: $G_j(q,p) = q_jp_j$ \vspace{-0.6cm}
	\item Focus-focus component: $ \begin{array}{l} \\ \begin{cases} 
	G_j(q,p) = q_j p_{j+1} - q_{j+1}p_j \\ G_{j+1}(q,p) =q_j p_j + q_{j+1} p_{j+1} 	
	\end{cases} \end{array}	$
%	$  \begin{cases} 
%	Q_j(q,p) = q_j p_{j+1} - q_{j+1}p_j \\ Q_{j+1}(q,p) =q_j p_j + q_{j+1} p_{j+1} 	
%	\end{cases}$
\end{itemize} 
\begin{de}
A semitoric system is a 4-dimensional integrable system $(M,\omega,F=(L,H))$ with two degrees of freedom such that $L$ is a proper map inducing an effective $\mbS^1$-Hamiltonian action on $M$ and $F$ only has non-degenerate singularities without hyperbolic components. If, moreover, for all $ \ell \in \R$, the fibre $L^{-1}(\ell)$ contains at most one focus-focus singularity, the system is called simple. 
\label{defst}
\end{de}
Note that semitoric systems can only have three types of singularities: elliptic-regular, elliptic-elliptic and focus-focus. 

\begin{ex}
Consider $M=\mbS^2 \times \R^2$ with Cartesian coordinates $(x,y,z,u,v)$ and symplectic form $\om=\rho_1\, \om_{\mbS^2} \oplus \rho_2 \, \om_{\R^2}$, where $\rho_1, \rho_2>0$ and $\om_{\mbS^2},\om_{\R^2}$ are the standard symplectic forms on the unit sphere and the plane respectively. The \emph{coupled spin-oscillator} is the integrable system $(M,\om,(L,H))$ defined by
$$ L(x,y,z,u,v) = \rho_1(z-1) + \dfrac{\rho_2}{2} (u^2+v^2),\qquad H(x,y,z,u,v) = \frac{1}{2} (xu+yv)$$ and it is a semitoric system.
\end{ex}

Two semitoric systems $(M_1,\om_1,F_1)$, $(M_2,\om_2,F_2)$ are said to be \emph{isomorphic} if there exists a pair $(\varphi,g)$, where $\varphi:(M_1,\om_1) \to (M_2,\om_2)$ is a symplectomorphism and $g:\R^2 \to \R$ is a smooth map such that $(L_2,H_2) = (L_1,g(L_1,H_1))$ and $\partial_{H_1} g >0$. 

\subsection{Symplectic invariants}
Pelayo \& \vungoc \cite{PV1, PV4} classified semitoric systems in terms of five symplectic invariants. This means, on the one hand, that if two semitoric systems have the same invariants, then there is an isomorphism between them and, on the other hand, that given an admissible list of invariants, a semitoric system with such invariants exists. We now introduce the invariants.

\subsubsection{Number of focus-focus points}
The first invariant is the number $\nff \in \N$ of singular points of focus-focus type. The coupled spin-oscillator has two singular points of rank 0, namely $(0,0,1,0,0)$ and $(0,0,-1,0,0)$. The first one is of focus-focus type and the second one is of elliptic-elliptic type, so $\nff=1$.
 
\subsubsection{Taylor series invariant}
\label{sstaylor}
The Taylor series invariant is a $\nff$-tuple of Taylor series $S:= (S^\infty_1,...,S^\infty_\nff)$, each of them associated to one of the focus-focus singularities $m_1,...,m_\nff \in M$ of our system. It is a semiglobal invariant, in the sense that it only depends on the dynamics around the focus-focus singular fibres. More precisely, fix ${r}\in \{1,..,\nff\}$ and consider an open saturated neighbourhood $ U_r \subseteq M$ of $m_r$. Using the Eliasson-Miranda-Zung normal form, we can find a pair $(\varphi_{r},\varrho_r)$, where $\varphi_{r}$ is a local symplectomorphism and $\varrho_r$ a local diffeomorphism, such that the following diagram commutes:
\begin{small}
\begin{center}
\begin{tikzpicture}
 \matrix (m) [matrix of math nodes,row sep=4em,column sep=6em,minimum width=2em]
 {
 U_{r} \subseteq M & F(U_{r}) \subseteq \R^2 \\
 V_{r} \subseteq \R^4 & G(V_{r}) \subseteq \R^2 \simeq \C \\};
 \path[-stealth]
 (m-1-1) edge node [left] {$\varphi_{r}$} 
 				node [below,rotate=90] {$\sim$} (m-2-1)
 			edge node [above] {$F$} (m-1-2)
 (m-1-2) edge node [left] {$\varrho_r$}
 				node [below,rotate=90] {$\sim$} (m-2-2)
 (m-2-1) edge node [above] {$G$} (m-2-2);
\end{tikzpicture}
\end{center}
\end{small}
We can choose $\varrho_r=(\varrho_{{r},1},\varrho_{{r},2})$ in such a way that $\varrho_{{r},1}(l,h)=l$ and $\partial_h \varrho_{{r},2}(l,h)>0$. We use $\varrho_r$ to extend this construction to the semiglobal neighbourhood $W_{r} = F^{-1}(F(U_r))$ by defining $\Phi_{r}:W_{r} \to G(V_{r})$ as $\Phi_{r} := \varrho_r \circ F$. 

If $z:= l+ij$ is a coordinate on $G(V_{r})$, where $j$ is the value of the function $G_2$ from the normal form, we can consider regular fibres $\Lam_{z} := \Phi_{r}^{-1}(z) \subset M$ for $z \neq 0$, which are tori. If we denote by $\ga_z$ the vanishing cycle of $\Lam_z$ and by $\delta_z$ the other cycle, then the actions of the system will be given by $\mcL_{r}(z) := \frac{1}{2\pi} \oint_{\ga_z} \varpi \equiv l$ and $\mcI_{r}(z) := \frac{1}{2\pi} \oint_{\delta_z} \varpi$, where $\varpi$ is a semiglobal primitive of the symplectic form. \vungoc \cite{Vu1} showed that
$$ 2\pi \mcI_{r}(z) = 2\pi \mcI_{r}(0) - \text{Im}(z \log z - z) + S_{r}(z),$$
where $S_{r}(z) = S_{r}(l,j)$ is a holomorphic function and therefore analytic. Let $S_{r}^\infty$ denote its Taylor series. By construction, it has no constant term and, by choosing a suitable determination of the complex logarithm, we can also impose $0 \leq \partial_l S_{r}(0) < 2\pi$. Repeating this construction for all focus-focus points gives us the invariant.

In practice, the Taylor series invariant is one of the most challenging to compute. Through some analytical considerations, the \emph{linear} terms of the invariant were computed for the coupled spin-oscillator (cf.\ Pelayo $\&$ \vungoc \cite{PV3}, Babelon \& Douçot \cite{BD}) and the coupled angular momenta for a \emph{specific} value of the parameters (cf.\ Le Floch \& Pelayo \cite{LFPe}). However, in order to compute the higher order terms and their dependence on the parameters of the system, a different approach using computers was needed. So far, higher order terms of the Taylor series invariant have been only computed for the coupled spin-oscillator (cf.\ Alonso \& Dullin \& Hohloch \cite{ADH}) and the coupled angular momenta (cf.\ Alonso \& Dullin \& Hohloch \cite{ADH2}). 

\subsubsection{Polygon invariant}
\label{sspolinv}
The polygon invariant allows us to compare the standard affine structure of $\R^2$ with the natural affine structure of $B=F(M)$ induced by the map $F$. 

For each ${r}=1,...,\nff$, let $c_{r} = (\lam_{r},\eta_{r}) = F(m_{r})$ be the critical value, $b_{\lam_{r}}$ the vertical line through $c_{r}$, $\eps_{r} \in \{-1,+1\}$ a choice of sign and $b_{\lam_{r}}^{\eps_{r}}$ the vertical half line that starts in $c_{r}$ and continues upwards if $\eps_{r}=+1$ or downwards otherwise. \vungoc \cite{Vu2} shows that there exists a homeomorphism onto its image $f:B \to \Delta \subseteq \R^2$ such that $\Delta$ is a rational polygon, $f|_{B \backslash \bigcup b_{\lam_{r}}^{\eps_{r}}}$ is a diffeomorphism into its image and $f$ preserves the first coordinate. 

This map, sometimes called the \emph{cartographic homeomorphism} (cf.\ Sepe \& \vungoc \cite{SV}), is not unique. It depends on the choice of signs $\eps = (\eps_1,...,\eps_\nff)$ and can be composed with a vertical translation or an affine transformation 
\begin{equation} T^k=\begin{pmatrix} 1 & 0 \\ k & 1  \end{pmatrix},\qquad k \in \Z.  \label{tmatrix} \end{equation}

This induces an action of the group $\{-1,+1\}^\nff \times \mcG$ on the space of rational polygons, where $\mcG = \{T^k | k \in \Z\}$. The polygon invariant is then defined as the orbit of the polygon $\Delta$ under this group action, together with the vertical lines meeting the critical points and the sign choice: 
$$(\{-1,+1\}^\nff \times \mcG) \cdot (\Delta,(b_{\lam_{r}})_{{r}=1}^\nff,(\eps_{r})_{{r}=1}^\nff).$$

\subsubsection{Height invariant}
The height invariant $\mathbf{h} = (h_1,...,h_\nff)$ is an $\nff$-tuple of positive real numbers. Let $f$ be a cartographic invariant as in \S \ref{sspolinv}.

We define the \emph{generalised toric momentum map} ${\mu}: M \to \Delta$ as the composition ${\mu} = f \circ F$. Observe that ${\mu}=({\mu}_1,{\mu}_2)$ satisfies ${\mu}_1(p) = L(p)$. The \emph{height} of the focus-focus singular point $m_{r}$ is then
$$ h_{r} = {\mu}_2(m_{r}) - \min_{p \in \Delta \cap b_{\lam_{r}}} \pi_2(p), $$
where $\pi_2:\R^2 \to \R$ is the projection onto the second coordinate. It can also be interpreted as the symplectic volume of the submanifold 
$$Y_{r}^- := \{ p \in | L(p)=L(m_{r}) \mbox{ and } H(p)< H(m_{r})\}.$$ 
The height invariant is independent of the choice of $f$.

\subsubsection{Twisting index invariant}
Roughly, to each polygon $\Delta$ obtained in \S \ref{sspolinv} we can associate an $\nff$-tuple of integers $(k_{r},...,k_\nff)$ by comparing the generalised toric momentum map $\mu$ to the local actions described in \S \ref{sstaylor}. More precisely, for each $i=1,..,\nff$ we can compare $\mu$ to the so-called \emph{privileged momentum map} $\nu_{r}=(\mcL_{r},\mcI_{r}) = (L,\mcI_{r})$ defined on the semiglobal neighbourhood $W_{r}$ and determine $k_{r}$ from the relation ${\mu} = T^{k_{r}} \circ \nu_{r}$, where $T^{k_{r}}$ is as in \eqref{tmatrix}. The indices $k_{r}$ are not independent of the choice of cartographic homeomorphism $f$. The group $\mcG$ acts by addition and the group $\{-1,+1\}^\nff$ does not act. The twisting-index invariant is conventionally written attached to the polygon invariant, so we consider in fact the orbit
$$(\{-1,+1\}^\nff \times \mcG) \cdot (\Delta,(b_{\lam_{r}})_{i=1}^\nff,(\eps_{r})_{r=1}^\nff,(k_{r})_{r=1}^\nff).\\[0.28cm]$$

\section{The coupled angular momenta: a paradigmatic example}

\subsection{Description of the system}
\label{ssdes}
The coupled angular momenta is one of the most important (families of) examples of semitoric systems. It describes the classical version of a non-trivial coupling of two quantum angular momenta. Sadovskii \& Zhilinskii \cite{SZ} observed in 1999 a redistribution of the energy levels in the quantum version of the system as they varied the coupling parameter. In particular, they showed that for certain values of the parameter, the system exhibits monodromy, an obstruction to finding global action-angle coordinates. 

Le Floch \& Pelayo \cite{LFPe} proved that the coupled angular momenta are semitoric for almost all values of the coupling parameter, constituting the first examples of semitoric systems with a compact phase space, in this case $\mbS^2 \times \mbS^2$. They partially computed the list of symplectic invariants, quantised the system and computed the associated joint spectrum in terms of Berezin-Toeplitz operators. Alonso \& Dullin \& Hohloch \cite{ADH2} generalised the results to all three parameters of the family and completed the computation of all symplectic invariants.

Let $R_1,R_2$ be positive real constants and consider the symplectic manifold $(M,\om)$  given by $M=\mbS^2 \times \mbS^2$ and $\om = -(R_1 \om_{\mbS^2} \oplus R_2  \om_{\mbS^2} )$, where $\om_{\mbS^2}$ is the standard symplectic form on the unit sphere. The \emph{coupled angular momenta} is the completely integrable system $(M,\om,(L,H))$ given by
\begin{equation}
\begin{cases}
L(x_1,y_1,z_1,x_2,y_2,z_2)\;:= R_1(z_1-1) + R_2 (z_2+1), \\
H(x_1,y_1,z_1,x_2,y_2,z_2):= (1-t) z_1 + t (x_1x_2 + y_1y_2 + z_1z_2) +2t-1,\\
\end{cases}
\label{cameq} 
\end{equation} where $t\in \R$ is the coupling parameter and $(x_i,y_i,z_i)$ are Cartesian coordinates on the unit sphere, $i=1,2$.

We will distinguish between the \emph{standard} case $R_1<R_2$, the \emph{Kepler} case $R_1=R_2$ and the \emph{reverse} case $R_1>R_2$. The discrete symmetries of \eqref{cameq} imply that the reverse case is isomorphic as a semitoric system to the standard case with the sign of $L$ inverted, which corresponds to a `time reversal' transformation (Prop.\ 3.2 of \cite{ADH2}). The Kepler case is isomorphic to the Kepler problem (Prop.\ 3.17 of \cite{ADH2}), that is the classical version of the hydrogen atom, if we express it in prolate spheroidal coordinates and perform symplectic reduction by the Hamiltonian flow. For a detailed discussion, see Dullin \& Waalkens \cite{DW}.

 Le Floch \& Pelayo \cite{LFPe} show that the system \eqref{cameq} is semitoric for all $t\in \R \backslash \{t^-,t^+\}$, where 
 \begin{equation}
 t^\pm := \dfrac{R_2}{2R_2+R_1 \mp 2\sqrt{R_1 R_2}}.
 \label{tpm}
 \end{equation} 
 
\begin{figure}[ht]
 \centering
 \includegraphics[width=14cm]{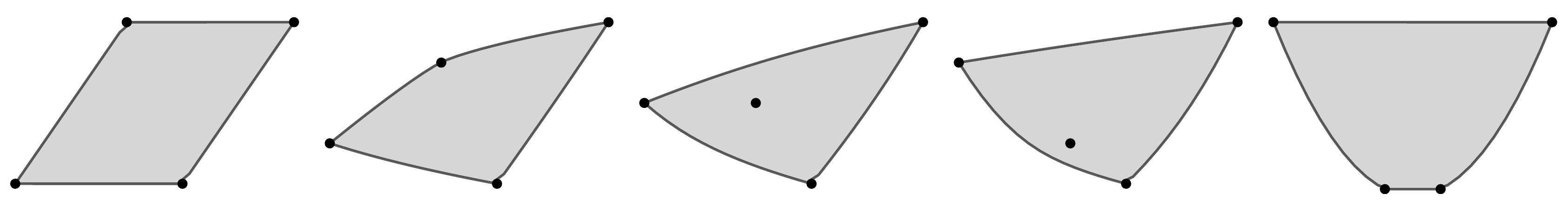}
 \caption{Image of the momentum map $F$ of the coupled angular momenta for $R_1=1$, $R_2=3/2$, from $t=0$ (left) to $t=1$ (right).}
 \label{vart}
\end{figure}

 The system has four singularities, namely $(z_1,z_2) = (\pm 1,\pm 1)$. The point $(1,-1)$ is of focus-focus type if $t^- < t < t^+$ and elliptic-elliptic if $t<t^-$ or $t>t^+$. The other three fixed points are always elliptic-elliptic. In Figure \ref{vart} we can see how, if $t$ changes, one of the elliptic-elliptic points --which are always in the border-- transforms into a focus-focus point --which are always interior-- and becomes again elliptic-elliptic by moving \emph{vertically}.

\begin{figure}[ht]
 \centering
 \includegraphics[width=15cm]{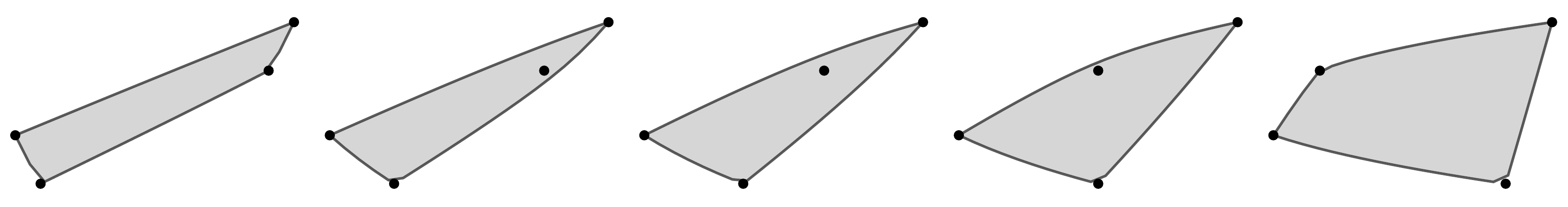}
 \caption{Image of the momentum map $F$ of the coupled angular momenta for $t=\frac{3}{10}$ and $R$ going from $R=\frac{1}{10}$ (left) to $R=5$ (right).}
 \label{varR}
\end{figure}

Traditionally, $t$ has been regarded as a parameter and $R:=R_2 / R_1$ as a scaling, but $R$ can also be understood as a parameter. If we invert \eqref{tpm}, we can obtain $R$ as a function of $t$. In Figure \ref{varR} we can see how, if $R$ changes, one of the elliptic-elliptic points  transforms into a focus-focus point and becomes again elliptic-elliptic by moving \emph{horizontally}. This figure was suggested to us by Holger Dullin.
 
\subsection{Symplectic invariants}
The coupled angular momenta is the only family of {\em compact} semitoric systems that has been completely classified. Summarising the recent progress by Alonso \& Dullin \& Hohloch \cite{ADH2}, the invariants are as follows.

\subsubsection{Number of focus-focus points}
From \S \ref{ssdes} it is clear that $\nff=0$ if $t<t^-$ or $t>t^+$ and $\nff=1$ if $t^-<t<t^+$. 
 
\subsubsection{Taylor series invariant}
The Taylor series symplectic invariant of the coupled angular momenta is given by
\begin{align*}
S(l,j) =& \, l \arctan \left( \dfrac{{R_2}^2 ( 2 t-1) - R_1 R_2 ( t+1) + {R_1}^2 t}{(R_1 - R_2)r_D} \right)+j \ln \left( \frac{4 {r_D}^3}{{R_1}^{1/2}{R_2}^{3/2}(1-t)t^2} \right) \\ \nonumber
+& \dfrac{l^2}{16 R_1R_2 {r_D}^3} \left( - {R_2}^4 ( 2 t-1)^3 -R_1 {R_2}^3 (32t^3-46 t^2+17 t-1)\right. \\ \nonumber &\hspace{2.2cm} \left. -3 {R_1}^2 {R_2}^2 t ( 4 t^2- 7 t+1)   
 -{R_1}^3 R_2 ( 5 t-3) t^2- {R_1}^4 t^3   \right) \\ \nonumber
+& \dfrac{lj}{8 R_1 R_2{r_D}^2} \left( (R_2 - R_1)({R_2}^2 (2t-1)^2   + 2 R_1 R_2 t ( 6 t-1)+{R_1}^2 t^2) \right)\\ \nonumber
+& \dfrac{j^2}{16 R_1 R_2 {r_D}^3} \left( 
  {R_2}^4 (  2 t-1)^3 
- R_1 {R_2}^3 (   16 t^3- 42 t^2+ 15 t+1)
\right. \\ \nonumber &\hspace{2.2cm} \left.
- {R_1}^2 {R_2}^2 t (28 t^2-3t-3) 
+ {R_1}^3 R_2 t^2 (  13 t-3)
+ {R_1}^4 t^3
 \right)
+ \mcO(3),
\end{align*} where 
\begin{small}
\begin{equation*}
 r_D: = \sqrt{-{R_2}^2 (1 - 2 t)^2 + 2 R_1 R_2 t - {R_1}^2 t^2}=\sqrt{({R_1}^2+4{R_2}^2)(t-t^-)(t^+-t)}.
\end{equation*} 
\end{small}
\subsubsection{Polygon and twisting index invariants}
\begin{figure}[ht!]
\centering
\begin{tabular}{m{1.5cm} m{3.6cm} m{3cm} m{3.6cm}}
\begin{center}$\mathbf{(k_1,\varepsilon_1)}$\end{center} &
\begin{center}$\mathbf{R_1 < R_2}$ \end{center} &
\begin{center}$\mathbf{R_1 = R_2}$ \end{center} &
\begin{center}$\mathbf{R_1 > R_2}$ \end{center} \\ 

\begin{minipage}{1.5cm} $k_1=-1$\\$\varepsilon_1=+1$ \end{minipage} &
\begin{center}\includegraphics[width=3.4cm]{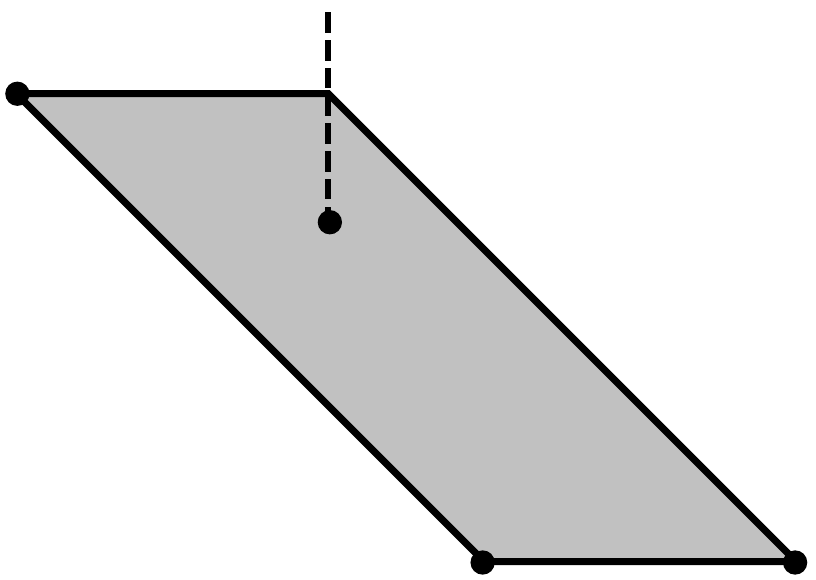}\end{center} &
\begin{center}\includegraphics[width=3.0cm]{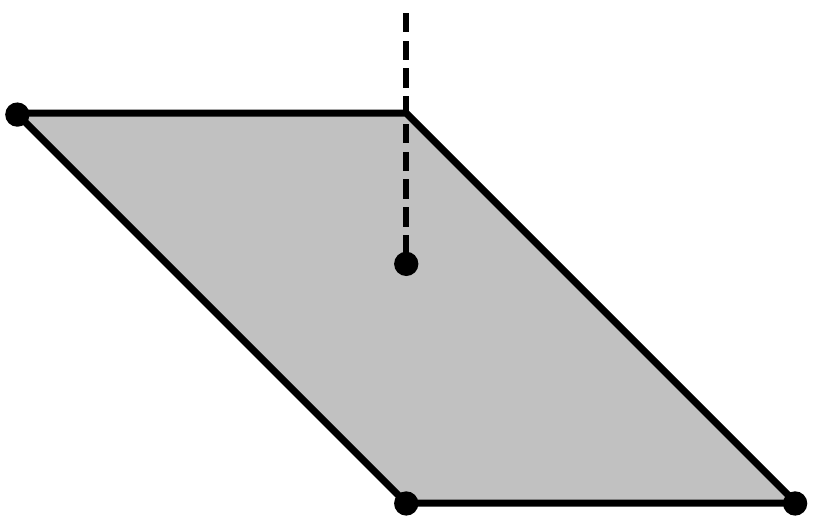}\end{center} &
\begin{center}\includegraphics[width=3.4cm]{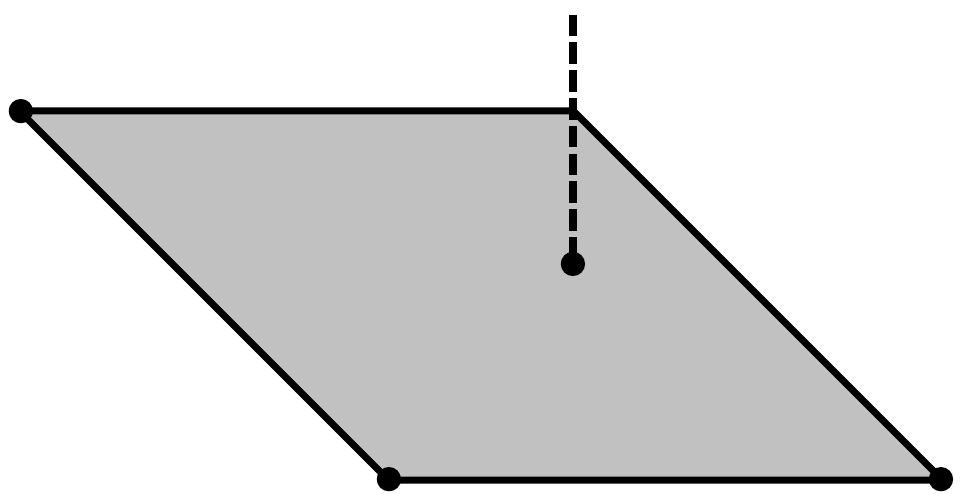}\end{center} \\

\begin{minipage}{1.5cm} $k_1=-1$\\$\varepsilon_1=-1$ \end{minipage} &
\begin{center}\includegraphics[height=1.8cm]{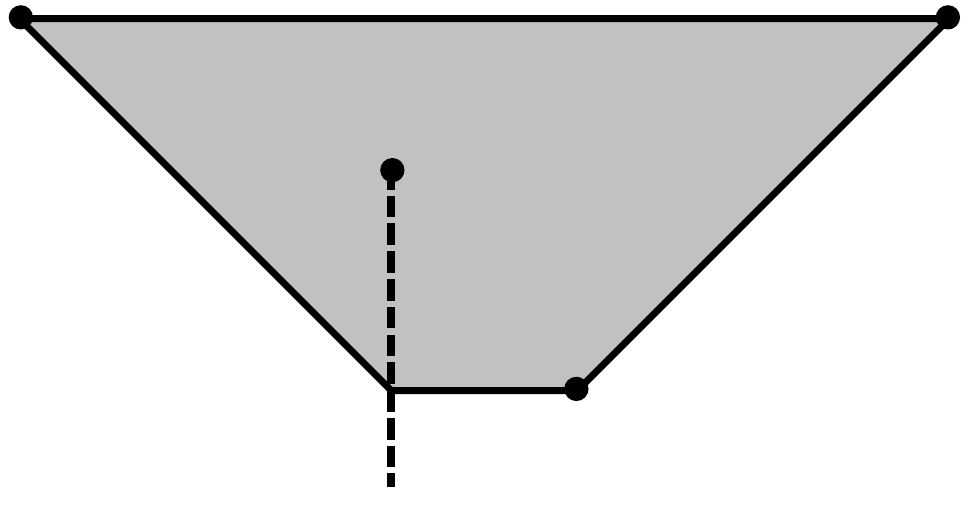}\end{center} &
\begin{center}\includegraphics[height=1.8cm]{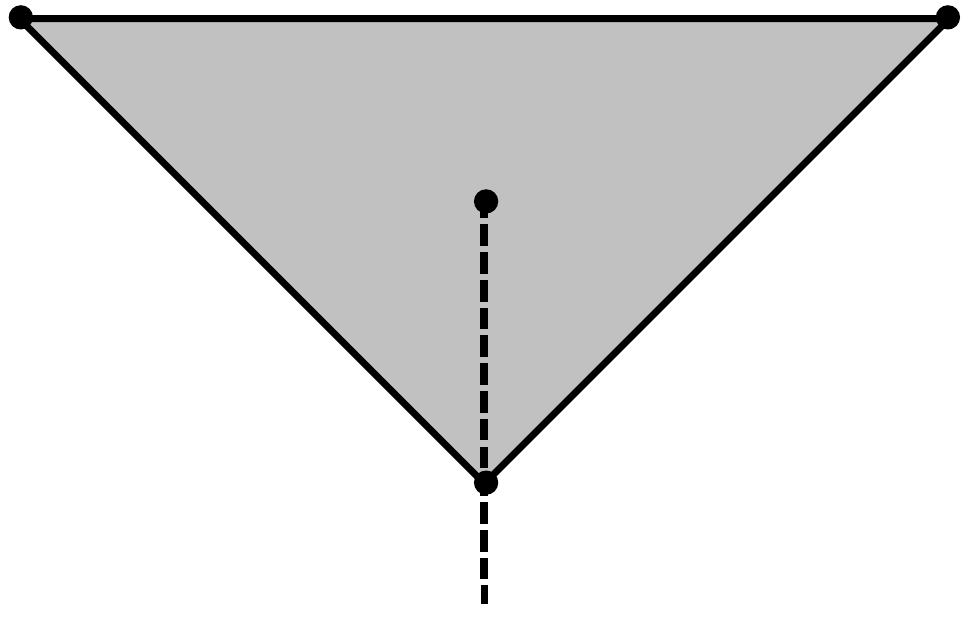}\end{center} &
\begin{center}\includegraphics[height=1.8cm]{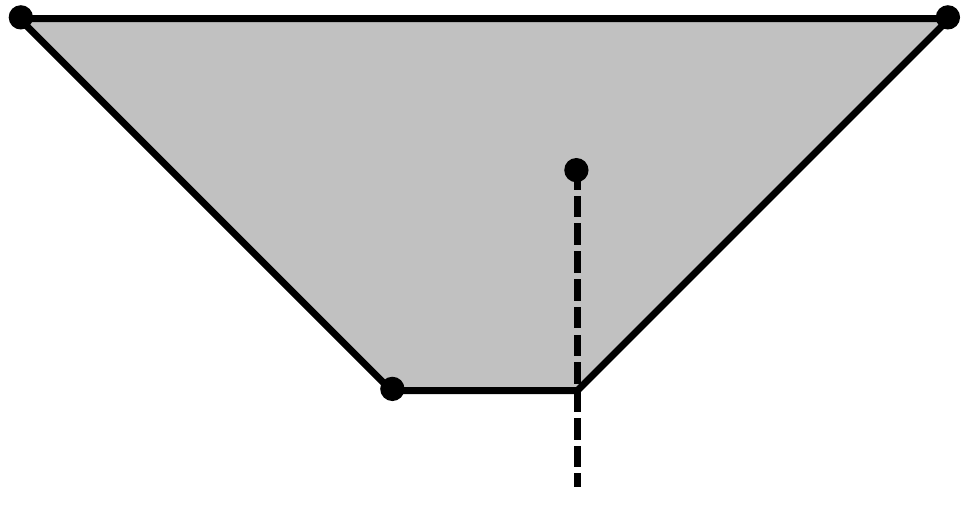}\end{center} \\

\begin{minipage}{1.5cm} $k_1=0$\\$\varepsilon_1=+1$ \end{minipage} &
\begin{center}\includegraphics[width=3.4cm]{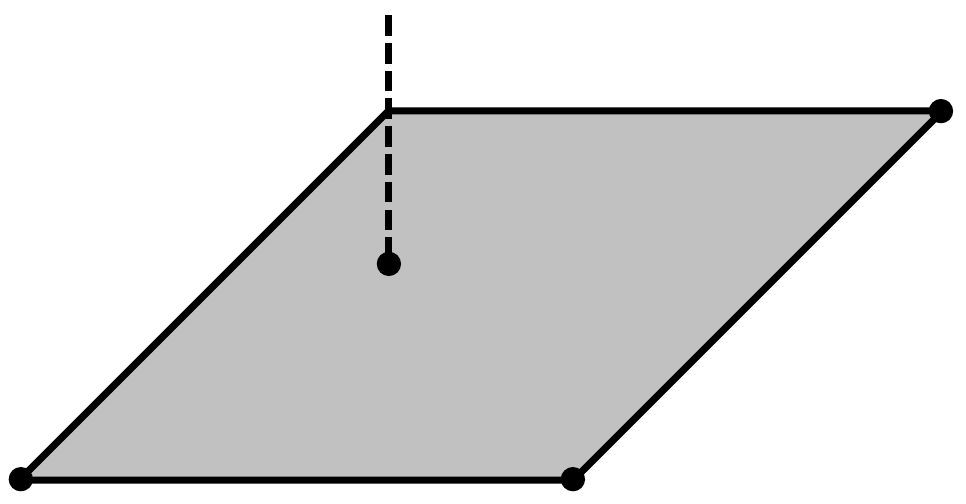}\end{center} &
\begin{center}\includegraphics[width=3.0cm]{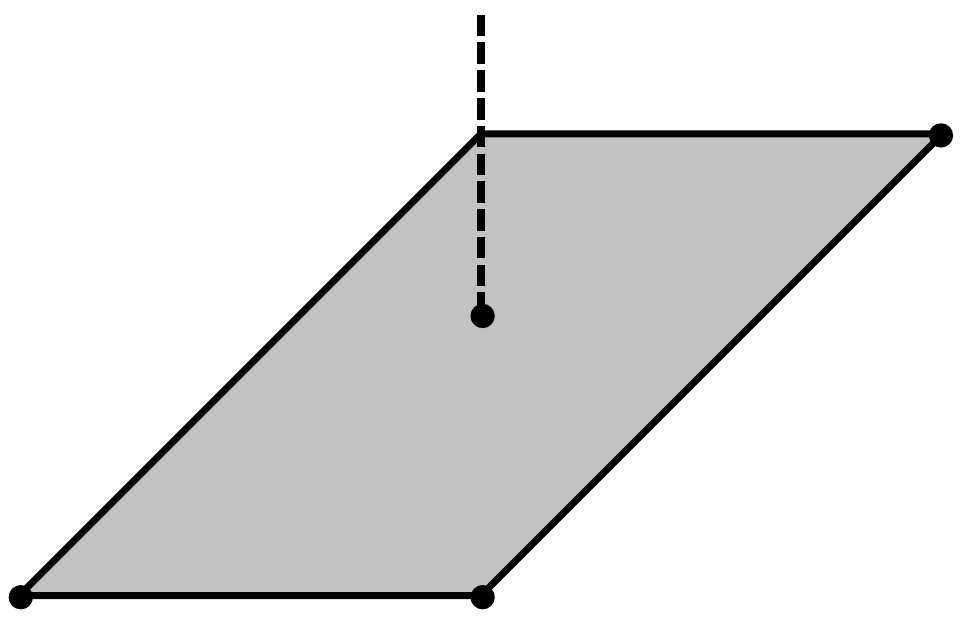}\end{center} &
\begin{center}\includegraphics[width=3.4cm]{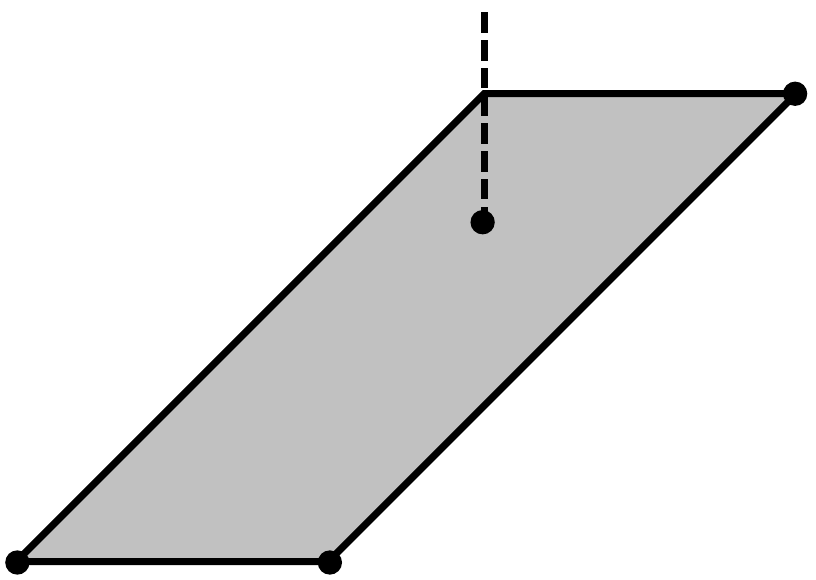}\end{center} \\
\end{tabular}
\caption{Some of the elements of the polygon invariant  for the cases $R_1<R_2$, $R_1=R_2$ and $R_1>R_2$ and $t^- < t < t^+$. The index $k_1$ is the associated twisting index to each polygon and $\varepsilon_1$ is the choice of sign.}
\label{poltwis}
\end{figure}

In Figure \ref{poltwis} we can see some of the polygons of the polygon invariant for the cases $R_1<R_2$, $R_1=R_2$ and $R_1>R_2$, together with their twisting index, for the case $t^- < t < t^+$. In this case the group acting is $\{-1,+1\} \times \mcG$. For $t<t^-$ we only have the polygons with $\eps_1=+1$ --without the vertical half-line indication-- and the group acting is $\mcG \simeq \Z$. The same happens for $t>t^+$, but in this case we only consider the polygons with $\eps_1=-1$.

\subsubsection{Height invariant}
The height invariant of the coupled angular momenta is
\begin{small}
\begin{align*}
h=2\min\{R_1,R_2\} + \dfrac{r_D}{\pi t} - \dfrac{2R_2}{\pi} \arctan \left( \dfrac{r_D}{R_2-tR_1} \right)-
\dfrac{2R_1}{\pi} \arctan \left( \dfrac{r_D}{R_2+tR_1-2R_2t} 
\right).
\end{align*}
\end{small}
\section{A new generation of examples}
In order to understand semitoric systems better, there was also a push for more general families of examples, such as containing more than one focus-focus points or with a compact underlying manifold other than $\mbS^2 \times \mbS^2$.

\subsection{A system with two focus-focus points}
\begin{figure}[ht]
 \centering
 \begin{tabular}{m{5cm} m{7cm}}
  \includegraphics[width=4.5cm]{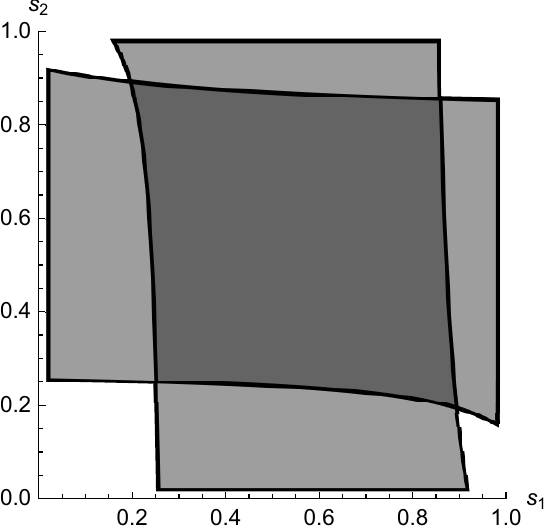} &
 \includegraphics[width=7cm]{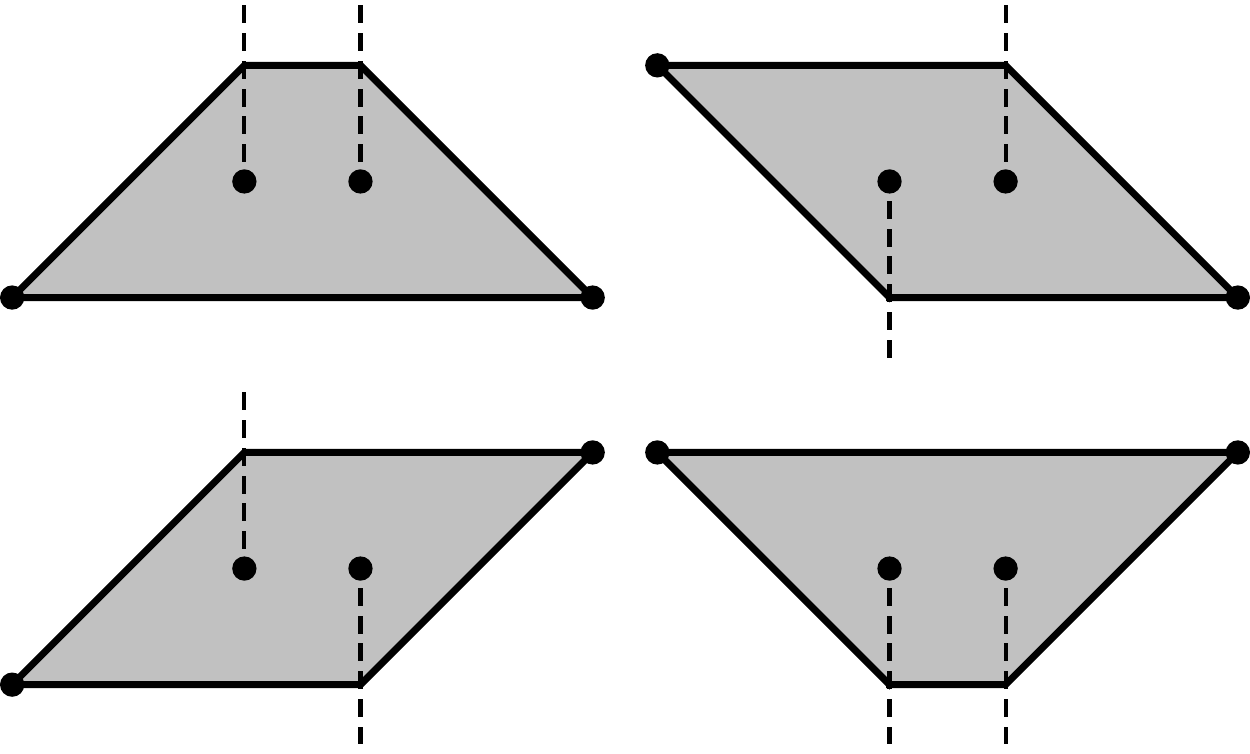}
 \end{tabular}
 \caption{Number of focus-focus (left) and some elements of the polygon invariant (right) of the system \eqref{2ffeq}.}
 \label{nff}
\end{figure}

The two classical families of examples coming from physics, the coupled spin-oscillator and the coupled angular momenta, have at most one singularity of focus-focus type. By introducing more parameters in the coupled angular momenta, Hohloch \& Palmer \cite{HP} managed to give the first example of semitoric system with two focus-focus singularities. More precisely, let $M=\mbS^2 \times \mbS^2$ and $\om = -(R_1 \om_{\mbS^2} \oplus R_2  \om_{\mbS^2} )$ be as before. Consider the system $(M,\om,(L,H))$ given by
\begin{equation}
\begin{cases}
L(x_1,y_1,z_1,x_2,y_2,z_2)\;:=& \hspace{-0.3cm} R_1 z_1  + R_2 z_2, \\
H(x_1,y_1,z_1,x_2,y_2,z_2):=& \hspace{-0.3cm} (1-s_1)(1-s_2)z_1 + s_1s_2 z_2 \\&+ s_1(1-s_2)(x_1x_2+y_1y_2+z_1z_2) \\&+ s_2(1-s_1)(x_1x_2+y_1y_2-z_1z_2),\\
\end{cases}
\label{2ffeq} 
\end{equation} where $0<R_1<R_2$ and $(s_1,s_2) \in [0,1]^2$. The system \eqref{2ffeq} is semitoric for all $(s_1,s_2)\in [0,1]^2\backslash \ga$, where $\ga$ is the union of certain four smooth curves. If $(s_1,s_2) \in \ga$, then at least one of the singularities is degenerate. These curves mark also the borders of the regions where the system has two (dark gray), one (light gray) or no (white) focus-focus singularity, as it can be observed in Figure \ref{nff}, left. 

Note that if we set $s_1=t$, $s_2=0$ in the system \eqref{2ffeq} then we recover the coupled angular momenta \eqref{cameq} up to a translation in the momentum map. Moreover, in the region of parameters $(s_1,s_2)$ for which the system has two focus-focus singularities, the system is unique in some sense. More precisely, Kane \& Palmer \& Pelayo \cite{KPP} showed that all semitoric systems on a compact manifold with two focus-focus points such that $L$ has isolated fixed points have the same polygon invariant, up to rescaling. Figure \ref{nff}, right, displays some of the polygons of this invariant.

\subsection{Semitoric families in Hirzebruch surfaces}
\label{sshirze}
An important breakthrough for generating examples of semitoric systems is the article by Le Floch \& Palmer \cite{LFPa}. Given a partial list of invariants, such as the number of focus-focus points and the polygon invariant, they looked for a strategy to find some simple examples of semitoric systems having those invariants. To this end, they introduced the concept of a \emph{semitoric family} and \emph{semitoric transition family}, which allowed them to find examples of semitoric systems in all Hirzebruch surfaces. This was the first time that examples were given on a compact manifold other than $\mbS^2 \times \mbS^2$. 

Hirzebruch surfaces can be defined in several equivalent ways, but here we are interested in the so-called Delzant construction (cf.\ Delzant \cite{Del} and Cannas da Silva \cite{Cds}). Given $\al, \be >0$ and $n \in \N$, the $n$-th Hirzebruch surface $(W_n(\al,\be),\om_{W_n(\al,\be)})$ is defined by taking $\C^4$ with the standard symplectic form $\om_{std}=\frac{i}{2} \sum_j du_j \wedge d\bar{u}_j$ and performing symplectic reduction at the level zero with respect to the Hamiltonian $\T^2$-action generated by the map
$$ N(u_1,u_2,u_3,u_4):= \dfrac{1}{2} \left( |u_1|^2 + |u_2|^2 + n|u_3|^2,|u_3|^2+|u_4|^2 \right) - \left(\al + \be n,\be \right).$$ 
All these Hirzebruch surfaces admit a toric system $(L_{std},H_{std})$ given by $L_{std}:= \frac{1}{2} |u_2|^2$ and $H_{std}:=\frac{1}{2}|u_3|^2$. Note also that, for $n=0$, we recover $W_0 (\al,\be) \simeq \mbS^2 \times \mbS^2$, the symplectic manifold underlying the coupled angular momenta. 

Let $(M,\om)$ be a four-dimensional symplectic manifold. Assume that $(M,\om,F_t)$ is a completely integrable system for all $0 \leq t \leq 1$, where $F_t := (L,H_t)$, $H_t:=H(t,\cdot)$ and $H:[0,1] \times M \to \R$ is a smooth map. Le Floch \& Palmer \cite{LFPa} call $(M,\om,F_t)$ a \emph{semitoric family} if there exist $k \in \N$ and $t_1, \ldots, t_k \in [0,1]$ such that $(M,\om,F_t)$ is semitoric for $t \in [0,1] \backslash \{t_1,\ldots,t_k\}$. Inspired by the coupled angular momenta, they define a so-called \emph{semitoric transition family} for a transition point $p \in M$ and transition times $t^-,t^+$ as follows. Consider a semitoric family $(M,\om,F_t)$ such that the point $p$ is a singularity of elliptic-elliptic type if $t<t^-$ or $t>t^+$ and of focus-focus type for $t^- < t < t^+$. Moreover, they require that if $p$ is a minimum (resp.\ maximum) of $H_0|_{L^{-1}(L(p))}$, then $p$ is a maximum (resp. minimum) of $H_1|_{L^{-1}(L(p))}$.

\begin{figure}[ht]
 \centering
 \includegraphics[width=6cm,trim=0 0 0 0]{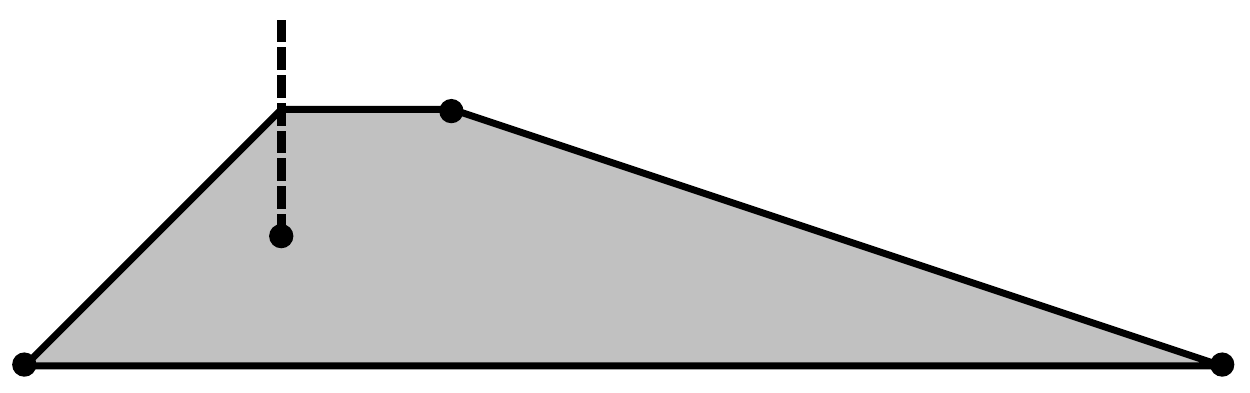} 
 \includegraphics[width=6cm,trim=0 0.8cm 0 0cm]{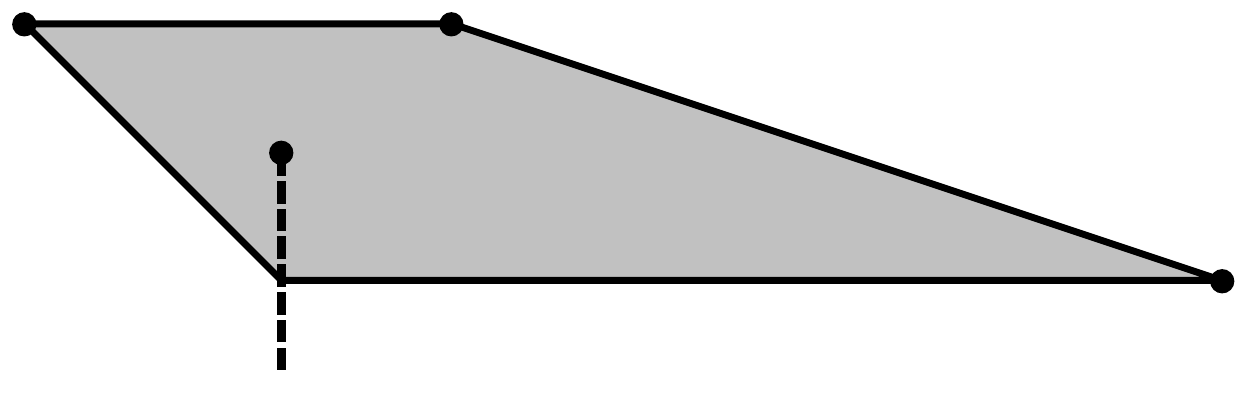}
 \caption{Some elements of the polygon invariant of $W_n(\al,\be)$. Left with $\eps_1=+1$ and right with $\eps_1=-1$.}
 \label{hirze}
\end{figure}

An important result (see Theorem 5.2 of \cite{LFPa}) is that, for all $ n \in \N$ and for all $ \al, \be >0$, there exists a semitoric transition family on $W_n(\al,\be)$ with degenerate points $0<t^-< t^+<1$. For $t<t^-$, the system has zero focus-focus points and its polygon invariant is the $\mcG$-orbit of the polygon in Figure \ref{hirze} left. The same happens for $t>t^+$ with the polygon on the right. For $t^- < t < t^+$, the system has one focus-focus point and the polygon invariant consists of the $\mcG$-orbits of both polygons in Figure \ref{hirze}. These systems are obtained by an alternating performance of $n$ blow-ups and blow-downs starting from the coupled angular momenta $W_0(\al',\be)$ for some $\al'>\al$.

Le Floch \& Palmer \cite{LFPa} also give some explicit examples of semitoric families in $W_1(\al,\be)$ and $W_2(\al,\be)$. Note that these examples are not the ones generated by the previous result. We reproduce two of them here: 
\begin{ex}[Thm.\ 6.2 of \cite{LFPa}]
Given $\al,\be>0$ and $0 < \ga < 1/(2 \sqrt{2\be})$, we define the system $F_t = (L_{std},H_t)$ on $W_1(\al,\be)$ by
$$ H_t := (1-2t) H_{std} + t \gamma X,\quad X:=\text{Re}(\bar{u}_1 u_3 \bar{u}_4).$$ This system is a semitoric transition family with at most one focus-focus point and transition times
$$ 0 <  \ t^- = \frac{1}{2(1+\ga \sqrt{2\be})} \ <  \ t^+ = \frac{1}{2(1-\ga \sqrt{2\be})} \ < 1.$$
\end{ex} 
\begin{ex}[Thm.\ 7.2 of \cite{LFPa}]
Let $\al,\be>0$ and define $\nu:=\be / \al$. Take also $0 < \ga < 1/(2\nu)$ and define $c:=2\ga \sqrt{\nu}$ so that $0<c<1$. Then the system in $W_2(\al,\be)$ defined by
$$L=\dfrac{1}{2}\left( |u_2|^2+|u_3|^2 \right),\quad H_t =(1-t)R + \dfrac{\be t (\ga X +(2L-\al-2\be)(R+\al+\be))}{\al(\al+2\be)},$$ where
$$ R := \dfrac{1}{2} \left( |u_3|^2-|u_4|^2 \right),\qquad X:=\text{Re}(\bar{u}_1\bar{u}_2u_3\bar{u}_4)$$ is a semitoric transition family with at most one focus-focus point and transition times
$$ t^- = \dfrac{1+2\nu}{1+(3+c)\nu},\quad t^+=\dfrac{1+2\nu}{1+(3-c)\nu}.$$
\end{ex}

\section{Advances beyond explicit examples and their invariants}

So far, this survey focused on the research progress in the study of examples of semitoric systems. However, there has also been significant advance in other aspects of the theory of semitoric systems and related research fields, as well as on the study of similar systems that do not satisfy some of the conditions of Definition \ref{defst}. In this section we summarise --far from giving an exhaustive list-- some of this progress.
 
Wacheux \cite{Wa} makes an extension of semitoric systems to dimension greater than four. Pelayo \& Ratiu \& \vungoc \cite{PRV, PRV2} generalise the theory of semitoric systems by relaxing the properness assumption while keeping connectivity of the fibres. Hohloch \& Sepe \& Sabatini \& Symington \cite{HSSS} work on so-called \emph{faithful semitoric systems} which are useful, in particular, under the point of view of categories of systems and surgeries. Convexity questions are studied by Ratiu \& Wacheux \& Zung \cite{RWZ}. Tentative steps towards hyperbolic singularities were made by Dullin \& Pelayo \cite{DP} and Le Floch \& Palmer \cite{LFPa}. Fans of semitoric systems in the sense of algebraic geometry and minimal models were studied by Kane \& Palmer \& Pelayo \cite{KPP, KPP2}.

Hohloch \& Sabatini \& Sepe \cite{HSS} show how the classification of semitoric systems relates to the classification of the inderlying $\mathbb S^1$-action in the sense of Karshon \cite{karshon}.
For progress in spectral theory and quantisation, we refer the reader to the lecture notes by Sepe \& \vungoc \cite{SV} and the articles Le Floch \& Pelayo \cite{LFPe}, Le Floch \& Pelayo \& \vungoc \cite{LFPV} and the references therein.

\bibliographystyle{alpha}
\bibliography{../bibtex/PhD}
\vspace{0.5cm}

\textbf{Jaume Alonso}\\
Department of Mathematics and Computer Science\\
University of Antwerp\\
Middelheimlaan 1\\
2020 Antwerp, Belgium\\
\textit{E-mail:} \texttt{jaume.alonsofernandez@uantwerpen.be}

\vspace{0.5cm}
\textbf{Sonja Hohloch}\\
Department of Mathematics and Computer Science\\
University of Antwerp\\
Middelheimlaan 1\\
2020 Antwerp, Belgium\\
\textit{E-mail:} \texttt{sonja.hohloch@uantwerpen.be}
\end{document}